\documentclass{amsart}

\newtheorem{theorem}{Theorem}
\newcommand{\bt}{\begin{theorem}}
\newcommand{\et}{\end{theorem}}

\newtheorem*{theoremNN}{Theorem}
\newcommand{\btNN}{\begin{theoremNN}}
\newcommand{\etNN}{\end{theoremNN}}

\newtheorem{lemma}{Lemma}
\newcommand{\bl}{\begin{lemma}}
\newcommand{\el}{\end{lemma}}

\newtheorem{corollary}{Corollary}
\newcommand{\bc}{\begin{corollary}}
\newcommand{\ec}{\end{corollary}}

\newtheorem{definition}{Definition}
\newcommand{\bdf}{\begin{definition}}
\newcommand{\edf}{\end{definition}}

\newtheorem{conjecture}{Conjecture}
\newcommand{\bconj}{\begin{conjecture}}
\newcommand{\econj}{\end{conjecture}}

\newtheorem*{conjectureNN}{Conjecture}
\newcommand{\bconjNN}{\begin{conjectureNN}}
\newcommand{\econjNN}{\end{conjectureNN}}

\newtheorem{example}{Example}
\newcommand{\bex}{\begin{example}}
\newcommand{\eex}{\end{example}}

\newtheorem{problem}{Problem}
\newcommand{\bprob}{\begin{problem}}
\newcommand{\eprob}{\end{problem}}

\newtheorem*{problemNN}{Problem}
\newcommand{\bprobNN}{\begin{problemNN}}
\newcommand{\eprobNN}{\end{problemNN}}

\newtheorem{oproblem}{Open Problem}
\newcommand{\boprob}{\begin{oproblem}}
\newcommand{\eoprob}{\end{oproblem}}

\newtheorem*{oproblemNN}{Open Problem}
\newcommand{\boprobNN}{\begin{oproblemNN}}
\newcommand{\eoprobNN}{\end{oproblemNN}}

\newcommand{\beq}{\begin{equation}}
\newcommand{\eeq}{\end{equation}}
\newcommand{\benum}{\begin{enumerate}}
\newcommand{\eenum}{\end{enumerate}}

\newcommand{\N}{\ensuremath{ \mathbf N }}

\newcommand{\Q}{\ensuremath{ \mathbf{Q} }}

\newcommand{\R}{\ensuremath{\mathbf R}}

\newcommand{\mbc}{\ensuremath{\mathbf c}}

\newcommand{\mbr}{\ensuremath{\mathbf r}}

\newcommand{\mcc}{\ensuremath{ \mathcal C}}

\newcommand{\mcr}{\ensuremath{ \mathcal R}}
\newcommand{\mcs}{\ensuremath{ \mathcal S}}

\newcommand{\bq}{\begin{eqnarray*}}
\newcommand{\eq}{\end{eqnarray*}}
\newcommand{\be}{\begin{eqnarray}}
\newcommand{\ee}{\end{eqnarray}}
\newcommand{\ba}{\begin{array}}
\newcommand{\ea}{\end{array}}
\newcommand{\bfr}{\begin{flushright}}
\newcommand{\efr}{\end{flushright}}

\newcommand{\bmat}{\left(\begin{matrix}}
\newcommand{\emat}{\end{matrix}\right)}
\newcommand{\bsmallmat}{\left(\begin{smallmatrix}}
\newcommand{\esmallmat}{\end{smallmatrix}\right)}

\DeclareMathOperator{\col}{\text{col}}

\DeclareMathOperator{\diag}{\text{diag}}

\DeclareMathOperator{\qand}{\quad\text{and}\quad}
\DeclareMathOperator{\qqand}{\qquad\text{and}\qquad}

\DeclareMathOperator{\row}{\text{row}}

\usepackage[all]{xy}\usepackage{amssymb,latexsym}

\title[Explicit Sinkhorn limits]{Matrix scaling, explicit Sinkhorn limits, and arithmetic}
\author{Melvyn B. Nathanson}
\address{Department of Mathematics\\Lehman College (CUNY)\\Bronx, NY 10468}
\email{melvyn.nathanson@lehman.cuny.edu}



\date{\today}

\begin{document}
\maketitle

\begin{abstract}
The process of alternately row scaling and column scaling a positive $n \times n$ matrix $A$ 
converges to a doubly stochastic positive $n \times n$ matrix $S(A)$, called the 
\emph{Sinkhorn limit} of $A$.  Exact formulae for the Sinkhorn limits of certain  
symmetric positive $3\times 3$ matrices are computed, and related problems in 
diophantine approximation are considered.  
\end{abstract}


\section{Doubly stochastic matrices and scaling}
Let $A = (a_{i,j})$ be an $m \times n$ matrix. 
For $i \in \{1,\ldots, m\}$, the $i$th \emph{row sum}\index{row sum} of $A$ is 
\[
\row_i(A) = \sum_{j=1}^n a_{i,j}.
\]
For $j \in \{1,\ldots, n\}$, the $j$th \emph{column sum}\index{column sum} of $A$ is 
\[
\col_j(A) = \sum_{i =1}^m a_{i,j}. 
\]
For example, the matrices
\[
\bmat
1 & -1 & 1 \\
-1 & 4 & -2 \\
1 & -2 & 2
\emat 
\qqand 
\bmat
2 & -5 & 4 \\
-9 & 7 & 3 \\
8 & -1 & -6
\emat
\]  
have row and column sums equal to 1.

An $n \times n$  matrix $ (u_{i,j})$ is \emph{diagonal} if $u_{i,j} = 0$ for all $i \neq j$.  
Let  $\diag(x_1,\ldots, x_n)$ denote the diagonal matrix whose $(i,i)$th coordinate 
is $x_i$ for  all $i \in \{1,2,\ldots, n\}$.  
The diagonal matrix $\diag(x_1,x_2,\ldots, x_n)$ 
is \emph{positive  diagonal} if $x_i > 0$ for all $i$.

The process of multiplying the rows of a matrix $A$ by scalars, 
or, equivalently, multiplying $A$ on the left by a diagonal matrix $X$, 
is called \emph{row-scaling}, 
and $X$ is called a \emph{row-scaling matrix}.  

The process of multiplying the columns of a matrix $A$ by scalars, 
or, equivalently, multiplying $A$ on the right by a diagonal matrix $Y$,  
is called \emph{column-scaling}, 
and $Y$ is called a \emph{column-scaling matrix}.

Let $A = (a_{i,j})$ be an $m \times n$ matrix.  
If $X = \diag(x_1,x_2,\ldots, x_m)$ and $Y = \diag(y_1,y_2,\ldots, y_n)$, then 
\[
XAY = 
\bmat
x_1a_{1,1} y_1 & x_1a_{1,2} y_2 & x_1a_{1,3} y_3 &  \cdots & x_1a_{1,n} y_n \\ 
x_2a_{2,1} y_1 & x_2a_{2,2} y_2 &x_2a_{2,3} y_3 &  \cdots & x_2a_{2,n} y_n \\ 
\vdots &&&& \vdots \\
x_m a_{m,1} y_1 & x_m a_{m,2} y_2 &x_m a_{m,3} y_3 &  \cdots & x_m a_{m,n} y_n 
\emat.
\]

The  $m \times n$ matrix  $ A = (a_{i,j})$ is \emph{positive} if $a_{i,j}>0$ for all $i$ and $j$, 
and \emph{nonnegative} if $a_{i,j} \geq 0$ for all $i$ and $j$. 
The matrix $A = (a_{i,j})$ is \emph{row stochastic}\index{row stochastic} 
if $A$ is nonnegative and  $\row_i(A) = 1$ for all $i \in \{1,\ldots, m\}$.  
The  matrix $A$ is \emph{column stochastic}\index{column stochastic} 
if $A$ is nonnegative and  $\col_j(A) = 1$ for all $j \in \{1,\ldots, n\}$.  
The matrix $A$ is \emph{doubly stochastic}\index{doubly stochastic} 
if it is both row and column stochastic.  
For example, 
the matrices
\[
\bmat 1/3 & 1/3 & 1/3 \\ 
1/3 & 1/3 & 1/3   \\ 
1/3 & 1/3 & 1/3 \emat 
\qqand 
\bmat 1/2 & 1/3 & 1/6 \\ 
1/6 & 1/2 & 1/3   \\ 
1/3 & 1/6 & 1/2 \emat,
\]  
are doubly stochastic.  

If the $m \times n$ matrix $A$ is doubly stochastic, then 
\[
m = \sum_{i=1}^m \row_i(A) = \sum_{i=1}^m \sum_{j=1}^n a_{i,j} 
=  \sum_{j=1}^n \sum_{i=1}^m a_{i,j} 
= \sum_{j=1}^n \col_j(A) = n
\]
and so $A$ is a square matrix.

Let $A = (a_{i,j})$ be an $m \times n$ matrix with positive row sums, 
that is, $\row_i(A) > 0$ for all $i \in \{1,\ldots, m\}$.  
Let $X(A) = \diag(1/\row_1(A), \ldots, 1/\row_m(A) )$ denote the $m\times m$ diagonal matrix 
whose $i$th diagonal coordinate is $1/\row_i(A)$, and let 
\[
\mcr(A) = X(A) A.
\]
We have 
\[
\mcr(A)_{i,j} = \frac{a_{i,j}}{\row_i(A)}
\]
and so
\[
\row_i(\mcr(A)) 
= \sum_{j=1}^n \mcr(A)_{i,j} 
= \sum_{j=1}^n\frac{a_{i,j}}{\row_i(A)} 
= \frac{\row_i(A)}{\row_i(A)} = 1
\]
for all $i \in \{1,\ldots, m\}$.  
Therefore, $\mcr(A)$ is a row stochastic matrix.

Similarly, let $Y(A) = \diag(1/\col_1(A), \ldots, 1/\col_n(A) )$ denote the $n\times n$ diagonal matrix 
whose $j$th diagonal coordinate is $1/\col_j(A)$, and let 
\[
\mcc(A) = A Y(A).
\]
We have 
\[
\mcc(A)_{i,j} = \frac{a_{i,j}}{\col_j (A)}
\]
and so
\[
\col_j(\mcc(A)) 
= \sum_{j=1}^n \mcc(A)_{i,j} = \sum_{i=1}^m \frac{a_{i,j}}{\col_j(A)} 
= \frac{\col_j(A)}{\col_j(A)} = 1
\]
for all $j \in \{1,\ldots, n\}$.  
Therefore, $\mcc(A)$ is a column stochastic matrix.

For example, if 
\[
A = \bmat 1 & 2 & 3 \\ 4 & 5 & 6 \emat
\]
then the matrix 
\begin{align*}
\mcr(A) = X(A) A 
& = \bmat 1/6 & 0 \\ 0 & 1/15 \emat  \bmat 1 & 2 & 3 \\ 4 & 5 & 6 \emat 
 =  \bmat 1/6 & 1/3 & 1/2 \\ 4/15 & 1/3 & 2/5 \emat 
\end{align*}
is row stochastic, and the matrix 
\begin{align*}
\mcc(A) = A Y(A) 
& = \bmat 1 & 2 & 3 \\ 4 & 5 & 6 \emat  \bmat 1/5 & 0 & 0\\ 0 &1/7 & 0 \\ 0 & 0 & 1/9 \emat 
 =  \bmat 1/5 & 2/7 & 1/3 \\ 4/5 & 5/7 & 2/3 \emat 
\end{align*}
is column stochastic.

In this paper we study doubly stochastic matrices.  

The following results (due to Sinkhorn~\cite{sink64}, Knopp-Sinkhorn~\cite{sink-knop67}, 
Menon~\cite{meno67},  Letac~\cite{leta74}, Tverberg~\cite{tver76}, and others) are classical.

\bt
Let $ A = (a_{i,j})$ be an $n \times n$ matrix with $a_{i,j}>0$ for all $i,j \in \{1,\ldots, n\}$.
\benum
\item[(i)]
There exist positive diagonal $n\times n$ matrices $ X$ and $ Y$ 
such that $ X A  Y$ is doubly stochastic.
\item
If $ X$, $ X'$, $ Y$, and $ Y'$ are  positive diagonal $n\times n$ matrices such that 
both  $ X A  Y$ and $ X' A  Y'$ are doubly stochastic, then 
$ X A  Y =  X' A  Y'$ and there exists $\lambda > 0$ 
such that $ X' = \lambda  X$ and $ Y' = \lambda^{-1}  Y$.  

The unique doubly stochastic matrix $XAY$ is called the \emph{Sinkhorn limit} of A,
and denoted $S(A)$.  

\item
Let A\ be a positive symmetric $n \times n$ matrix.  There exists a unique positive diagonal matrix X\ 
such that $X A X$ is doubly stochastic. 
\eenum
\et

\bt
Let $\mcs_n$ be the set of positive doubly stochastic matrices.
Let $\R^{n}_{>0}$ (resp. $\R^{n-1}_{>0}$) be the set of positive $n$-dimensional 
(resp. $(n-1)$-dimensional)  vectors.
Consider 
\[
\Omega = \R^{n}_{>0} \times \mcs_n \times \R^{n-1}_{>0}
\]
 as a subset of $\R^{n^2+2n-1}$ with the subspace topology.  
Consider the set $M_n^+$  of positive $n\times n$ matrices 
as a subset of $\R^{n^2}$ with the subspace topology.
The function from $\Omega$  to $M_n^+$ 
defined by 
\[
\bmat
x_1 \\ \vdots \\ x_{n-1} \\ x_n 
\emat , S, 
\bmat
y_1 \\ \vdots \\ y_{n-1} \\ 1 
\emat 
\mapsto \diag(x_1,\ldots, x_{n-1}, x_n) \  S \   \diag(y_1,\ldots, y_{n-1}, 1)
\]
is a homeomorphism.  
\et

\bt
Let A\ be a positive $n\times n$ matrix.  
Construct  sequences of positive matrices  $(A_{\ell})_{\ell=0}^{\infty}$ 
and $(A'_{\ell})_{\ell=0}^{\infty}$ 
and sequences of positive diagonal matrices $(X_{\ell})_{\ell=0}^{\infty}$ 
and $(Y_{\ell})_{\ell=0}^{\infty}$ 
as follows:
Let 
\[
A_0 = A.
\]
Given the matrix $A_{\ell}$, let 
\[
X_{\ell} = X(A_{\ell}) =  \diag\left( \frac{1}{\row_1(A_{\ell})},  
\frac{1}{\row_2(A_{\ell})},  \ldots,  \frac{1}{\row_n(A_{\ell})}\right)
\]
be the row-scaling matrix of $A_{\ell}$, and let 
\[
A'_{\ell} = X_{\ell} A_{\ell}.
\]
The matrix $A'_{\ell}$ is row stochastic.  
Let 
\[
Y_{\ell} =  Y(A'_{\ell}) 
= \diag\left( \frac{1}{\col_1(A)},  \frac{1}{\col_2(A)},  \ldots,  \frac{1}{\col_n (A)}\right)
\]
be the column-scaling matrix of $A'_{\ell}$, and let 
\[
A_{\ell+1} =   A'_{\ell}Y_{\ell}.
\]
The matrix $A_{\ell +1}$ is column stochastic.  
There exist positive diagonal matrices X\ and Y\ such that 
\[
\lim_{\ell\rightarrow \infty} X_{\ell} = X, \qquad
\lim_{\ell\rightarrow \infty} Y_{\ell} = Y
\]
and the $n \times n $ matrix 
\[
S(A) = XAY = \lim_{\ell\rightarrow \infty} A_{\ell} = \lim_{\ell\rightarrow \infty} A'_{\ell}
\]
is doubly stochastic.  
\et

This process of obtaining a doubly stochastic matrix $S(A)$  from a positive matrix $A$ 
by row and column scaling is called \emph{alternate minimization}.  

It is an open problem to compute explicitly the Sinkhorn limit of a positive $n\times n$ matrix.  
This is known  for $2\times 2$ matrices (Nathanson~\cite{nath18aaa}).  
In this paper we compute explicit Sinkhorn limits for certain symmetric $3 \times 3$ matrices,
and discuss connections with diophantine approximation.


\section{Experimental data}       \label{Sinkhorn:section:calculation}
Here are some computational results.    
Using Maple, we row scale and then column scale the matrix, iterate this process 20 times,  
and print the resulting matrix.

\[
\bmat 
2 & 1 & 1 \\
1 & 1 & 1  \\
1 & 1 & 1 \emat
\rightarrow 
\bmat
0.4384471874 & 0.2807764064 & 0.2807764064 \\
0.2807764064 & 0.3596117968 & 0.3596117968 \\ 
0.2807764064 & 0.3596117968 & 0.3596117968 
\emat
\]

\[
\bmat 
1 & 1 & 1 \\
1 & 2 & 2 \\
1 & 2 & 2 \emat
\rightarrow 
\bmat
0.4384471873 & 0.2807764065 & 0.2807764065 \\
0.2807764064 & 0.3596117968 & 0.3596117968 \\
0.2807764064 & 0.3596117968 & 0.3596117968 
\emat
\]

\[
\bmat
2 & 1 & 1 \\
1 & 2 & 1 \\
1 & 1 & 1
\emat
\rightarrow
\bmat 
0.4648162417 & 0.2324081208 & 0.3027756377 \\ 
0.2324081208 & 0.4648162417 & 0.3027756377 \\ 
0.3027756380 & 0.3027756380 & 0.3944487245 
\emat
\]

\[
\bmat
2 & 2 & 1 \\
2 & 1 & 1 \\
1 & 1 & 1
\emat
\rightarrow 
\bmat
0.3274800021  & 0.4125989480 & 0.2599210499 \\ 
0.4125989480 & 0.2599210499 & 0.3274800021 \\
0.2599210499 & 0.3274800021 & 0.4125989480 
\emat
\]

\[
\bmat
2 & 2 & 1 \\
2 & 1 & 1 \\
1 & 1 & 2
\emat = 
\bmat
0.3451802671 & 0.4435474272 & 0.2112723057 \\ 
0.4435474272 & 0.2849733008 & 0.2714792720 \\ 
0.2112723057 & 0.2714792720 & 0.5172484223 
\emat.
\]
In these calculations, the alternate minimization algorithm generates approximately 
doubly stochastic matrices of four different shapes: 
\[
\bmat
a & b & b \\
b & c & c  \\
b & c & c  
\emat,
\qquad
\bmat
a & b & c \\
b & a & c \\
c & c & d
\emat,
\qquad
\bmat
a & b & c \\
b & c & a \\
c & a & b
\emat,
\qquad
\bmat
a & b & c \\
b & d & e \\
c & e & f 
\emat. 
\]

\section{Permutation matrices} 
Let $S_n$ be the group of permutations of the set $\{1,2,\ldots, n\}$.    
For every $\sigma \in S_n$, define the $n\times n$ 
\emph{permutation matrix} $P_{\sigma}$ 
as follows:
\beq       \label{Sinkhorn:PermutationMatrix}
\left(P_{\sigma}\right)_{i,j}  
 = 
\begin{cases}
1 & \text{if $j = \sigma(i)$} \\
0 & \text{ if $j \neq \sigma(i)$.}
\end{cases} 
\eeq
Equivalently, 
\[
\left(P_{\sigma^{-1}}\right)_{i,j}  
 = 
\begin{cases}
1 & \text{if $i = \sigma(j)$} \\
0 & \text{ if $i \neq \sigma(j)$.}
\end{cases}
\]    
Thus, \[
\left(P_{\sigma}\right)_{i,j}  = \delta_{\sigma(i),j} = \delta_{i,\sigma^{-1}(j)}
\]
where  $\delta_{i,j}$ is the Kronecker delta.  
The $i$th row of $P_{\sigma}$ is row $\sigma(i)$ 
of the $n\times n$ identity matrix $I_n$, and 
the $j$th column of $P_{\sigma^{-1}}$ is column $\sigma(j)$ of $I_n$.

For every $n\times n$ matrix $A$, 
the $i$th row of the matrix $P_{\sigma}A$  
is row $\sigma(i)$ of $A$, and  the $j$th column of the matrix $AP_{\sigma^{-1}}$ 
is column $\sigma(j)$ of $A$.
Thus,  $P_{\sigma}A$ is a matrix 
constructed from $A$ by the $\sigma$-permutation of the rows of $A$, and  
$AP_{\sigma^{-1}}$ is a matrix 
constructed from $A$ by the $\sigma$-permutation of the columns of $A$. 

For example, if $\sigma = (1,2,3)$, then 
\[
P_{\sigma}A = \bmat 
0 & 1 & 0 \\
0 & 0 & 1 \\
1 & 0 & 0 
\emat 
\bmat 1 & 2 & 3 \\ 4 & 5 & 6 \\ 7 & 8 & 9 \emat
= 
\bmat 
4 & 5 & 6  \\  7 & 8 & 9 \\ 1 & 2 & 3 
\emat 
\]
and 
\[
A P_{\sigma^{-1}} = 
\bmat 1 & 2 & 3 \\ 4 & 5 & 6 \\ 7 & 8 & 9 \emat 
\bmat 
0 & 0 & 1 \\
1 & 0 & 0 \\
0 & 1 & 0 
\emat 
= 
\bmat 
 2 & 3 & 1 \\ 5 & 6 & 4 \\ 8 & 9 & 7 
\emat.
\]

\bl
For all permutations $\sigma,\tau \in S_n$,  
\beq       \label{Sinkhorn:PermutationMatrixCommute}
P_{\sigma}P_{\tau} = P_{\tau\sigma}
\eeq
and
\beq       \label{Sinkhorn:PermutationMatrixTranspose}
P_{\sigma}^t = P_{\sigma^{-1}}.
\eeq
\el

\begin{proof}
Let $i,j \in \{1,2,\ldots, n\}$.  
Applying~\eqref{Sinkhorn:PermutationMatrix} with $j=k$, we obtain 
\begin{align*}
\left(P_{\sigma}P_{\tau}\right)_{i,j} 
& = \sum_{k=1}^n \left(P_{\sigma}\right)_{i,k} \left(P_{\tau}\right)_{k,j} \\
& = \sum_{k=1}^n \delta_{\sigma(i),k} \left(P_{\tau}\right)_{k,j} \\
& =  \left(P_{\tau}\right)_{\sigma(i),j} \\
& = \delta_{\tau\sigma(i),j} \\
& = \left(P_{\tau\sigma}\right)_{i,j}.
\end{align*}
This proves~\eqref{Sinkhorn:PermutationMatrixCommute}.

For the transpose of $P_{\sigma}$, we have 
\begin{align*}
\left( P_{\sigma}^t \right)_{i,j} 
& = \left(P_{\sigma} \right)_{i,j} = p_{j,i} =
\begin{cases}
1 & \text{if $i = \sigma(j)$ } \\
0 & \text{ if $i \neq \sigma(j)$ }
\end{cases} \\
& = \left(P_{\sigma^{-1}} \right)_{i,j}.
\end{align*}
This proves~\eqref{Sinkhorn:PermutationMatrixTranspose}. 
\end{proof}

For example, if $\sigma = (1,2,3)$ and $\tau = (1,2)$, then $ \tau\sigma = (2,3)$.  
We have 
\[
P_{\sigma} = 
\bmat 
0 & 1 & 0 \\
0 & 0 & 1  \\
1 & 0 & 0 
\emat
\qqand
P_{\tau} = 
\bmat 
0 & 1 & 0  \\
1 & 0 & 0 \\
0 & 0 & 1 
\emat
\]
and
\[
P_{\sigma} P_{\tau} = 
\bmat 
0 & 1 & 0  \\ 
0 & 0 & 1  \\
1 & 0 & 0 
\emat
 \bmat 
0 & 1 & 0  \\
1 & 0 & 0 \\
0 & 0 & 1 
\emat 
 = \bmat 
1 & 0 & 0  \\
0 & 0 & 1 \\
0 & 1 &  0
\emat
= P_{\sigma\tau}.
\]

For $k, \ell \in \{1,2,\ldots, m\}$ with $k \neq \ell$, let 
$\tau \in S_m$ be the transposition defined by 
\[
\tau(k) = \ell, \qquad
\tau(\ell) = k
\]
and
\[
\tau(i) = i \qquad \text{
for all $i \neq k,\ell$.}
\]

Let $A = (a_{i,j})$ be an $m \times n$ matrix.  
The $m \times m$ permutation matrix $P_{\tau}$ interchanges 
rows $k$ and $\ell$ of $A$, as follows:
For all $i \in \{1,\ldots, m \}$  and $j \in \{1,\ldots, n\}$, 
\begin{align*}
\left(P_{\tau}A\right)_{k,j} & = a_{\ell,j} \\
\left(P_{\tau}A\right)_{\ell,j} & = a_{k,j} \\
\left(P_{\tau}A\right)_{i,j} & = a_{i,j} \qquad \text{ if $i \neq k,\ell$.}
\end{align*}
It follows that 
\begin{align*}
\mcr(P_{\tau}A)_{i,j} 
& = \frac{ \left(P_{\tau}A\right)_{i,j} }{ \row_i\left(P_{\tau}A\right) } = 
\begin{cases}
a_{\ell,j} / \row_{\ell}(A) & \text{ if $i =k$}\\
a_{k,j} /\row_k(A)  & \text{ if $i = \ell$} \\
a_{i,j}  /\row_i(A)& \text{ if $i \neq k,\ell$}
\end{cases} \\
& = \left( P_{\tau} \mcr(A)\right)_{i,j} 
\end{align*}
and so 
\beq       \label{Sinkhorn:row-transpose}
\mcr(P_{\tau}A) = P_{\tau} \mcr(A).  
\eeq

Let $\sigma$ be a permutation in $S_m$, and let $P_{\sigma}$ 
be the corresponding $m \times m$ permutation matrix.  
Every permutation $\sigma \in S_m$ is a product of transpositions, and so 
there is a sequence of transpositions $\tau_1,\ldots, \tau_{q-1}, \tau_q$ such that 
\[
\sigma = \tau_1\cdots \tau_{q-1} \tau_q
\]
and 
\[
P_{\sigma} = P_{\tau_1} \cdots P_{\tau_{q-1}} P_{\tau_q}.
\]
Applying identity~\eqref{Sinkhorn:row-transpose} recursively, we obtain 
\begin{align*}
\mcr(P_{\sigma} A ) 
& = \mcr(P_{\tau_1} P_{\tau_2}  \cdots P_{\tau_{q-1}} P_{\tau_q}A ) \\ 
& = P_{\tau_1} \mcr( P_{\tau_2}  \cdots P_{\tau_{q-1}} P_{\tau_q}A ) \\
& = \cdots \\
& = P_{\tau_1}P_{\tau_2}  \cdots P_{\tau_{q-1}}  \mcr( P_{\tau_q}A ) \\
& = P_{\tau_1}P_{\tau_2}  \cdots P_{\tau_{q-1}} P_{\tau_q} \mcr( A ) \\
& = P_{\sigma} \mcr( A ).
\end{align*}
This proves that, for all permutations $\sigma \in  S_m$, 
\beq       \label{Sinkhorn:R-row-permute}
\mcr(P_{\sigma}A) = P_{\sigma} \mcr(A)  
\eeq

Similarly, 
\beq       \label{Sinkhorn:R-column-permute}
\mcr(A Q_{\sigma}) = \mcr(A)  Q_{\sigma}
\eeq

\beq       \label{Sinkhorn:C-row-permute}
\mcc(P_{\sigma}A) = P_{\sigma} \mcc(A)  
\eeq

\beq       \label{Sinkhorn:C-column-permute}
\mcc(A Q_{\sigma}) = \mcc(A)  Q_{\sigma}
\eeq

For example, let 
\[
A = \bmat 1 & 2 & 3 \\ 4 & 5 & 6 \\ 7 & 8 & 9 \emat.
\]
Consider the permutation $\sigma = (3,2,1) \in S_3$ and its associated 
permutation matrix 
\[
P_{\sigma} = \bmat 0 & 0 & 1 \\ 1 & 0 & 0 \\ 0 & 1 & 0 \emat.
\]
We have 
\begin{align*}
\mcr(P_{\sigma} A ) 
& = \mcr\left( P_{\sigma} \bmat 1 & 2 & 3 \\ 4 & 5 & 6 \\ 7 & 8 & 9 \emat \right) 
= \mcr  \bmat  7 & 8 & 9  \\ 1 & 2 & 3 \\ 4 & 5 & 6 \emat  \\
& = \bmat  7/24 & 8/24 & 9/24  \\ 1/6 & 2/6 & 3/6 \\ 4/15 & 5/15 & 6/15 \emat
\end{align*}
and 
\begin{align*}
P_{\sigma} \mcr( A ) 
& = 
P_{\sigma} \mcr  \bmat 1 & 2 & 3 \\ 4 & 5 & 6 \\ 7 & 8 & 9 \emat 
= P_{\sigma} \bmat   1/6 & 2/6 & 3/6 \\ 4/15 & 5/15 & 6/15 \\ 7/24 & 8/24 & 9/24  \emat \\
& = \bmat  7/24 & 8/24 & 9/24  \\ 1/6 & 2/6 & 3/6 \\ 4/15 & 5/15 & 6/15 \emat.
\end{align*}

\bt
Let $A$ be an $m \times n$ matrix.  If $P$ and $Q$ are permutation matrices, then 
\begin{align*}
\mcr(PA) PA  & = P \mcr(A)A \\
\mcr(AQ) AQ & = \mcr(A) A Q \\
PA \mcc(PA) & = P A \mcc(A) \\ 
AQ \mcc(AQ) & = A \mcc(A) Q.
\end{align*}
\et

\begin{proof}
It suffices to prove this for transpositions.  

Interchanging two rows of a matrix and row scaling is the same as row scaling and then 
interchanging the rows.  

Interchanging two rows of a matrix and column scaling is the same as column scaling and then 
interchanging the rows.  

Interchanging two columns of a matrix and row scaling is the same as row scaling and then 
interchanging the columns.  

Interchanging two columns of a matrix and column scaling is the same as column scaling and then 
interchanging the columns.  
\end{proof}


\bt                         \label{Sinkhorn:theorem:PQ}
Let $A$ be an $n \times n$ positive matrix.  For all permutation matrices $P$ and $Q$, 
\[
S(PAQ) = PS(A)Q.
\]
\et

\begin{proof}
Let $\left( A^{(\ell)} \right)_{\ell=0}^{\infty}$ be the alternate minimization sequence of matrices 
constructed from $A = A^{(0)}$.  
For all $\ell \geq 0$, we have
\[
A^{(2\ell+1)} = \mcc \left( A^{(2\ell)}  \right) 
\]
\[
A^{(2\ell+2)} = \mcr \left( A^{(2\ell+1)}  \right) 
\]
and
\[
\lim_{\ell \rightarrow \infty} A^{(\ell)} = S(A).
\]
For every permutation matrix $P$, we have 
\[
(PA)^{(1)} = \mcc(PA) = P\mcc(A) = P A^{(1)}
\]
\[
\left( PA \right)^{(2)} = \mcr\left( (PA)^{(1)}  \right)  
= \mcr\left( PA^{(1)}  \right)  
= P\mcr\left( A^{(1)}  \right)  = P A^{(2)}
\]
\[
\left( PA \right)^{(3)} = \mcc\left( (PA)^{(2)}  \right)  
= \mcc\left( PA^{(2)}  \right)  
= P\mcc\left( A^{(2)}  \right)  = P A^{(3)}.
\]
Continuing inductively, we obtain 
\[
\left( PA \right)^{(\ell)}   = P A^{(\ell)}
\]
for all $\ell \in \N_0$, and so 
\begin{align*}
S(PA) & = \lim_{\ell \rightarrow \infty} \left( PA \right)^{(\ell)} 
\lim_{\ell \rightarrow \infty} PA^{(\ell)} \\
& = P \lim_{\ell \rightarrow \infty} A^{(\ell)} = PS(A).
\end{align*}
Similarly, for every permutation matrix $Q$, we have 
\[begin{align*}]
S(QA)  = S(A)Q.
\]
Therefore, 
\[
S(PAQ) = PS(AQ) = PS(A)Q.
\]
This completes the proof.  
\end{proof}

\bt              \label{Sinkhorn:theorem:transpose}
For every positive $n\times n$ matrix $A$, 
\[
S\left( A^t \right) = S(A)^t.
\]
\et
\begin{proof}
Let $X$ and $Y$ be diagonal matrices such that 
\[
S(A) = XAY.
\]
We have $X^t = X$, $Y^t = Y$, and 
\[
S(A)^t = \left(XAY\right)^t = Y^t A^t X^t = Y A^t X.  
\]
If $S(A)$ is doubly stochastic, then $S(A)^t$ is doubly stochastic.  
The uniqueness theorem implies that 
\[
S\left( A^t \right)= Y A^t X = S(A)^t.
\]
This completes the proof.  
\end{proof}

\bt              \label{Sinkhorn:theorem:dilate}
Let $\lambda > 0$.  
For every positive $n\times n$ matrix $A$, 
\[
\mcc (\lambda A) = \mcc(A), \qquad \mcr(\lambda A) = \mcr(A), 
\]
and
\[
S\left( \lambda A\right) = S(A).
\]
\et
\begin{proof}
Klar.
\end{proof}

Here is an example of permutation and dilation equivalence.  
Let 
\[
A = \bmat 2 & 2 & 2 \\ 3 & 2 & 2 \\2 & 2 & 3 \emat.
\]
Dilating $A$ by  $\lambda = 1/2$, we obtain
\[
\lambda A = \bmat 1 & 1 & 1 \\ 3/2 & 1 & 1 \\1 & 1 & 3/2 \emat.
\]
Multiplying by the permutation matrices 
\[
P = \bmat 0  & 0 & 1 \\ 0 & 1 & 0  \\ 1 & 0 & 0 \emat
\qqand
Q = \bmat   0 & 1 & 0 \\ 0 & 0 & 1  \\1 & 0 & 0 \emat
\]
we obtain 
\[
B = P(\lambda A) Q 
=  \bmat 3/2 & 1 & 1 \\ 1 & 3/2 & 1 \\ 1 & 1 & 1 \emat 
=  \bmat K & 1 & 1 \\ 1 & K & 1 \\1 & 1 & 1 \emat
\]
with $K = 3/2$.  
Equivalently, 
\[
A = \lambda^{-1}P^{-1} B Q^{-1} 
\]
and
\[
S(A) =  \lambda^{-1} P^{-1} S(B) Q^{-1}.
\]
Thus, the Sinkhorn limit of $B$ determines the Sinkhorn limit of A.

\section{The $MBN$ matrix}
Let $k$, $\ell$, and $n$ be positive integers such that $k+\ell = n$.  
Let $M$, $B$, and $N$ be positive real numbers.  
Consider the $n \times n$ symmetric matrix 
\beq            \label{Sinkhorn:MBN}
A = \bmat
M & M & \cdots & M & B & B & \cdots & B \\
M & M & \cdots & M & B & B & \cdots & B \\
\vdots &&& \vdots  & \vdots  &&& \vdots \\ 
M & M & \cdots & M & B & B & \cdots & B \\
 B & B & \cdots & B & N & N & \cdots & N \\
  B & B & \cdots & B & N & N & \cdots & N \\
\vdots &&& \vdots  & \vdots  &&& \vdots \\ 
 B & B & \cdots & B & N & N & \cdots & N \\
\emat
\eeq
in which the first $k$  rows are equal to 
\[
(\underbrace{M,M,\ldots, M}_{k},\underbrace{ B, B, \ldots, B}_{\ell})
\]
and the last $\ell$ rows  are equal to 
\[
(\underbrace{B,B,\ldots, B}_{k}, \underbrace{ N,N, \ldots, N}_{\ell}).
\]
Let $X  = \diag(x_1,x_2,x_3, \ldots, x_n ) $  be the unique positive $n \times n$  diagonal matrix 
such that the alternate minimization limit  $S(A) = X  A X $ is doubly stochastic. 
Thus, the matrix 
\[
S(A) 
 = \bmat
M x_1^2  & Mx_1x_2 & \cdots & M x_1x_k & B x_1x_{k+1} & B x_1x_{k+2} & \cdots & B x_1x_{n}\\
M x_2x_1  & Mx_2^2 & \cdots & M x_2x_k & B x_2x_{k+1} & B x_2x_{k+2} & \cdots & B x_2x_{n}\\
\vdots &&& \vdots  & \vdots  &&& \vdots \\ 
M x_k x_1  & Mx_kx_2 & \cdots & M x_k^2 & B x_kx_{k+1} & B x_kx_{k+2} & \cdots & B x_kx_{n}\\
B x_{k+1} x_1  & B x_{k+1} x_2 & \cdots & B x_{k+1}x_k & N x_{k+1}^2& N x_{k+1}x_{k+2} & \cdots & N x_{k+1}x_{n}\\
B x_{k+2} x_1  & B x_{k+2} x_2 & \cdots & B x_{k+2}x_k & N x_{k+2}x_{k+1}& N x_{k+2}^2 & \cdots & N x_{k+2}x_{n}\\
\vdots &&& \vdots  & \vdots  &&& \vdots \\ 
B x_n x_1  & B x_n x_2 & \cdots & B x_n x_k & N  x_nx_{k+1}& N x_nx_{k+2} & \cdots & N x_{n}^2
\emat
\]
satisfies 
\[
x_i\left( M \sum_{j=1}^k x_j + B\sum_{j=k+1}^n x_j \right) = 1 
\qquad \text{for $i = 1,2,\ldots k$} 
\]
and
\[
x_i\left( B \sum_{j=1}^k x_j + N\sum_{j=k+1}^n x_j \right) = 1 
\qquad \text{for $i = k+1, k+2, ,\ldots k+\ell = n$.} 
\]
It follows that $x_i = x_1$ for  $i = 1,2,\ldots k$ and  $x_i = x_n$ for  $i = k+1, k+2,\ldots n$.  
Let $x_1 = x$ and $x_n = y$.  Define the diagonal matrix 
\[
X = \diag(\underbrace{x,x,\ldots, x}_{k}, \underbrace{ y,y, \ldots, y}_{\ell}).
\] 
We obtain 
\begin{align}                
S(A)     \label{Sinkhorn:MBN-limit}  
& =  \bmat
M x^2& M x^2& \cdots & Mx^2& B  xy & B xy & \cdots & B xy \\
M x^2& M x^2& \cdots & M x_1^2 & B xy & B xy& \cdots & B xy \\
\vdots &&& \vdots  & \vdots  &&& \vdots \\ 
M x^2 & Mx^2& \cdots & M x_1^2 & B  xy & B  xy & \cdots & B xy \\
B  xy & B  xy& \cdots & B  xy & N y^2  & N  y^2 & \cdots & N  y^2 \\
B  xy  & B  xy& \cdots & B  xy & N y^2 & N  y^2  & \cdots & N  y^2 \\
\vdots &&& \vdots  & \vdots  &&& \vdots \\ 
B  xy & B  xy & \cdots & B xy & N   y^2 & N  y^2  & \cdots & N x_{n}^2
\emat                           \\
& =
\bmat
a  & a & \cdots & a & b & b & \cdots & b \\
a  & a & \cdots & a & b & b & \cdots & b \\
\vdots &&& \vdots  & \vdots  &&& \vdots \\ 
a  & a & \cdots & a & b & b & \cdots & b \\
 b & b & \cdots & b & c & c & \cdots & c \\
 b & b & \cdots & b &  c & c & \cdots & c \\
\vdots &&& \vdots  & \vdots  &&& \vdots \\ 
 b & b & \cdots & b & c & c & \cdots & c              \nonumber
 \emat                       
\end {align}
where 
\begin{align}       \label{Sinkhorn:MBN-a}
a & = Mx^2                            \\ 
b & = Bxy  = \frac{1-ka}{\ell}               \label{Sinkhorn:MBN-b}     \\ 
c & = Ny^2  = \frac{1-kb}{\ell} =  \frac{\ell-k + k^2 a}{\ell^2}.          \label{Sinkhorn:MBN-c}    
\end{align}
Because $S(A)$ is row stochastic, we have 
\beq            \label{Sinkhorn:MBN-1}
x \left(kMx+ \ell B y \right) = 1 
\eeq
and
\beq            \label{Sinkhorn:MBN-2}
y \left( kBx + \ell N y \right) = 1. 
\eeq
Equation~\eqref{Sinkhorn:MBN-1} gives 
\[
y =\frac{1}{\ell B} \left( \frac{1}{x} -kMx \right).
\]
Inserting this into equation~\eqref{Sinkhorn:MBN-2} and rearranging gives  
\beq            \label{Sinkhorn:MBN-3}
k^2 M \left( MN - B^2 \right) x^4 - \left( nB^2 + 2k (MN - B^2) \right) x^2+ N = 0 
\eeq 

If $MN-B^2 = 0$, then 
\[
x^2 = \frac{N}{nB^2} = \frac{1}{nM} 
\]
and $M x^2 = a = b = c = 1/n$.  Thus, $S(A)$ is the $n \times n$ doubly stochastic matrix 
with every coordinate equal to $1/n$.

If $MN - B^2 \neq 0$, then~\eqref{Sinkhorn:MBN-3} is a quadratic equation in $x^2$. 
We obtain 
\begin {align*}
x^2 
& = \frac{ 2k (MN - B^2)  + nB^2 \pm B\sqrt{ 4k\ell (MN-B^2) + n^2B^2}}{ 2k^2 M (MN-B^2)} \\
& =  \frac{ 1}{kM} + \frac{ nB^2 \pm B\sqrt{ 4k\ell MN  + (k-\ell)^2 B^2}}{ 2k^2 M (MN-B^2)} \\
\end {align*}
and
\begin {align*}
a  = M x^2 
& =  \frac{ 2k (MN - B^2) + nB^2 \pm B\sqrt{ 4k\ell (MN-B^2) + n^2B^2}}{ 2k^2(MN-B^2)} \\
& =  \frac{ 1}{k} + \frac{ nB^2 \pm B\sqrt{ 4k\ell MN  + (k-\ell)^2 B^2}}{ 2k^2 (MN-B^2)} \\
& =  \frac{ 1}{k} + \frac{ n \pm  \sqrt{ 4k\ell MN/B^2  + (k-\ell)^2 }}{ 2k^2  (MN/B^2 -1)}.  
\end {align*}

Recall that $ka+\ell b = 1$ and so $a< 1/k$.
If $MN > B^2$, then 
\[
n  = k+\ell <  \sqrt{ 4k\ell MN/B^2  + (k-\ell)^2 }.
\]
If $MN < B^2$, then 
\[
n  = k+\ell >  \sqrt{ 4k\ell MN/B^2  + (k-\ell)^2 }.
\]
In both cases, we obtain
\[                         
a  =  \frac{ 1}{k} + \frac{ n - \sqrt{ 4k\ell MN/B^2  + (k-\ell)^2 }}{ 2k^2  (MN/B^2 -1)}.  
\]
We obtain $b$ from~\eqref{Sinkhorn:MBN-b}  and $c$ from~\eqref{Sinkhorn:MBN-c}.  

\bt         \label{Sinkhorn:theorem:MBN-1-1}
The Sinkhorn limit of the $MBN$ matrix~\eqref{Sinkhorn:MBN} is the 
doubly stochastic matrix $S(A)$ defined by~\eqref{Sinkhorn:MBN-limit}.
The  matrix $S(A)$ depends only on the ratio $MN/B^2$.
\et

\begin{proof}
This follows immediately from~\eqref{Sinkhorn:MBN-a},~\eqref{Sinkhorn:MBN-b}, 
and~\eqref{Sinkhorn:MBN-c}.
\end{proof}

For example, the matrices 
\[
\bmat 
2 & 5 & 5 \\
5 & 3 & 3 \\
5 & 3 & 3
\emat,
\qquad
\bmat 
6 & 5 & 5 \\
5 & 1 & 1 \\
5 & 1 & 1
\emat,
\qqand
\bmat 
6/25 & 1 & 1 \\
1 & 1 & 1 \\
1 & 1 & 1
\emat,
\]
have the same Sinkhorn limit with $a = -37/38 + 5\sqrt{73}/38$.

Theorem~\ref{Sinkhorn:theorem:MBN-1-1} explains why, 
in Section~  \ref{Sinkhorn:section:calculation}, the matrices 
$
\bmat 
2 & 1 & 1 \\
1 & 1 & 1  \\
1 & 1 & 1 \emat
$
and
$\bmat 
1 & 1 & 1 \\
1 & 2 & 2 \\
1 & 2 & 2 \emat$
have the same Sinkhorn limits.

Let $\left( A^{(r)} \right)_{r=1}^{\infty}$ be a sequence of $MBN$ matrices 
such that $\lim_{r\rightarrow \infty} MN/B^2 = \infty$.
Let 
\[
S\left( A^{(r)} \right) = 
\bmat
a^{(r)}  & a^{(r)} & \cdots & a^{(r)} & b^{(r)} & b^{(r)} & \cdots & b^{(r)} \\
a^{(r)}  & a^{(r)} & \cdots & a^{(r)} & b^{(r)} & b^{(r)} & \cdots & b^{(r)} \\
\vdots &&& \vdots  & \vdots  &&& \vdots \\ 
a^{(r)}  & a^{(r)} & \cdots & a^{(r)} & b^{(r)} & b^{(r)} & \cdots & b^{(r)} \\
 b^{(r)} & b^{(r)} & \cdots & b^{(r)} & c^{(r)} & c^{(r)} & \cdots & c^{(r)} \\
 b^{(r)} & b^{(r)} & \cdots & b^{(r)} &  c^{(r)} & c^{(r)} & \cdots & c^{(r)} \\
\vdots &&& \vdots  & \vdots  &&& \vdots \\ 
 b^{(r)} & b^{(r)} & \cdots & b^{(r)} & c^{(r)} & c^{(r)} & \cdots & c^{(r)} \\
 \emat.
\]
We have 
\[
\lim_{r \rightarrow \infty} a^{(r)} = \frac{1}{k}, \qquad
\lim_{r \rightarrow \infty} b^{(r)}  = 0,   
\qquad 
\lim_{r \rightarrow \infty} c^{(r)} = \frac{1}{\ell}
\]
and
\[
\lim_{r\rightarrow \infty} S\left( A^{(r)} \right) = 
\bmat
1/k  & 1/k & \cdots & 1/k & 0 & 0 & \cdots & 0 \\
1/k  & 1/k & \cdots & 1/k & 0 & 0 & \cdots & 0 \\
\vdots &&& \vdots  & \vdots  &&& \vdots \\ 
1/k  & 1/k & \cdots & 1/k & 0 & 0 & \cdots & 0 \\
 0 & 0 & \cdots & 0 & 1/\ell & 1/\ell & \cdots & 1/\ell \\
 0 & 0 & \cdots & 0 &  1/\ell & 1/\ell & \cdots & 1/\ell \\
\vdots &&& \vdots  & \vdots  &&& \vdots \\ 
 0 & 0 & \cdots & 0 & 1/\ell & 1/\ell & \cdots & 1/\ell 
 \emat.
\]
Similarly, let $\left( A^{(r)} \right)_{r=1}^{\infty}$ be a sequence of $MBN$ matrices 
such that $\lim_{r\rightarrow \infty} MN/B^2 = 0$.  
It follows from~\eqref{Sinkhorn:MBN-a}  that 
\[
\lim_{r \rightarrow \infty} a^{(r)} 
=  \frac{ 1}{k} - \frac{ k+\ell  - |k-\ell | }{ 2k^2},
\]
If $k \leq \ell$, then 
\[
\lim_{r \rightarrow \infty} a^{(r)} = 0, 
\qquad 
\lim_{r \rightarrow \infty} b^{(r)} = \frac{1}{\ell}, 
\qquad
\lim_{r \rightarrow \infty} c^{(r)} =  \frac{\ell-k}{\ell^2}.  
\]
If $k > \ell$ , then 
\[
\lim_{r \rightarrow \infty} a^{(r)} =  \frac{k - \ell}{k^2},  
\qquad 
\lim_{r \rightarrow \infty} b^{(r)} =  \frac{k^2 - k + \ell}{2k^2}, 
\qquad
\lim_{r \rightarrow \infty} c^{(r)} =  \frac{k^2 + k - \ell}{2k^2}.
\]


\section{$3 \times 3$  symmetric matrices  and 
their doubly stochastic shapes}

Let $A$ and $B$ be $n\times n$ positive matrices.
We write $A \sim B$ if there exist $n\times n$  permutation matrices $P$ and $Q$ 
and $\lambda > 0$ such that 
\[
B = \lambda PAQ.
\]
It is straightforward to check that this is an equivalence relation.  
If $A\sim B$, then 
\[
S(B) = \lambda P S(A) Q.
\]
Thus, it suffices to compute the Sinkhorn limit of only one matrix in an equivalence class.

The goal is to compute the Sinkhorn limit of every $3\times 3$ symmetric positive  matrix
whose set of coordinates consists of two distinct  real numbers.
  
Let A\ be such a matrix with coordinates $a$ and $b$. 
There are 9 coordinate positions in the matrix, and so exactly one 
of the numbers $a$ and $b$ occurs at least five times.  
Suppose that the coordinate $a$ occurs five or more times.  Let $\lambda = 1/a$ 
and $K = b/a$.
The matrix $\lambda A$ has two distinct positive  coordinates $1$ and $K$, 
and $K$ occurs at most four times.  
There are seven equivalence classes of such matrices with respect to permutations and dilations.  
Here is the list, and, for each matrix, the shape of its Sinkhorn limit.
Note that $K$ is a positive real number and $K \neq 1$.

\benum

\item
\[
A_1 = \bmat
K & 1 & 1 \\
1 & K & 1 \\
1 & 1 & K
\emat
\qquad
S(A_1) = \bmat
a & b & b \\
b & a & b \\
b & b & a  
\emat
\]

\item
\[
A_2 = \bmat
K & 1 & 1 \\
1 & 1 & 1 \\
1 & 1 & 1
\emat
\qquad
S(A_2) = 
\bmat
a & b & b \\
b & c & c \\
b & c & c 
\emat
\]

\item
\[
A_3 = \bmat 
1 & 1 & 1 \\
1 & K & K  \\
1 & K & K  
\emat
\qquad
S(A_3) = \bmat
a & b & b \\
b & c & c \\
b & c & c
\emat
\]

\item
\[
A_4 = \bmat
1 & K & K \\
K & 1 & 1 \\
K & 1 & 1
\emat
\qquad
S(A_4) = \bmat
a & b & b \\
b & c & c \\
b & c & c
\emat
\]

\item
\[
A_5 = \bmat
K & 1 & 1 \\
1 & K & 1 \\
1 & 1 & 1
\emat
\qquad
S(A_5) = \bmat
a & b & c \\
b & a & c \\
c & c & d 
\emat
\]

\item
\[
A_6 = \bmat
K & K & 1 \\
K & 1 & 1 \\
1 & 1 & 1
\emat
\qquad
S(A_6) = \bmat
a & b &c \\
b & c & a \\
c & a & b 
\emat
\]

\item
\[
A_7 = \bmat
K  & K & 1 \\
K & 1 & 1  \\
1 & 1 & K
\emat
\qquad
S(A_7) = \bmat
a & b & c \\
b & d & e \\
c & e & f
\emat
\]

\eenum

\section{The matrix $A_1$}
The matrix 
\[
A_1 = \bmat
K & 1 & 1 \\
1 & K & 1 \\
1 & 1 & K
\emat
\]
is the simplest.  Just one row scaling or one column scaling produces 
the doubly stochastic matrix 
\[
A_1 \rightarrow S(A_1) = \bmat
K/(K+2)  & 1/(K+2) & 1/(K+2) \\
1/(K+2) & K/(K+2) & 1/(K+2) \\
1/(K+2) & 1/(K+2) & K/(K+2)  
\emat
\]
We have  $S(A_1) = X A_1 X$, where 
\[
X = \diag( \sqrt{1/(K+2)}, \sqrt{1/(K+2)}, \sqrt{1/(K+2)} ).
\]
Moreover,  
\[
\lim_{K\rightarrow \infty} S(A_1) = 
\bmat 
1 & 0 & 0 \\
0 & 1 & 0 \\
0 & 0 & 1 
\emat.
\]

\section{The matrices $A_2$, $A_3$,  and $A_4$}
These  are $MBN$ matrices.  
The matrix 
\[
A_2 = \bmat K & 1 & 1 \\ 1 & 1 & 1 \\ 1 & 1 & 1 \emat 
\] 
is an $MBN$ matrix with $k=1$, $\ell = 2$, $M=K$, and $B = N = 1$.

The matrix 
\[
A_3 = \bmat 
1 & 1 & 1 \\
1 & K & K  \\
1 & K & K  
\emat
\]
is an $MBN$ matrix with $k=1$, $\ell = 2$, $M=B=1$, and $N = K$.
Both matrices satisfy $MN/B^2 = K \neq 1$, and so they have the same Sinkhorn limit
\[
\bmat 
a & b & b \\
b & c & c \\
b & c & c
\emat
\]
with 
\begin{align}
a & =  \frac{2K+1 - \sqrt{8K  + 1}}{2(K-1)}     \label{Sinkhorn:A2-a}            \\
b& = \frac{ -3 + \sqrt{8K  + 1}}{4(K-1)}    \label{Sinkhorn:A2-b}            \\
c & = \frac{ 4K-1 - \sqrt{8K  + 1}}{8(K-1)}.   \label{Sinkhorn:A2-c}            
\end{align}

For example, if  $K=2$, then 
\[
A_2 = \bmat 2 & 1 & 1 \\ 1 & 1 & 1 \\ 1 & 1 & 1 \emat
\]
and 
\[
A_3 = \bmat 
1 & 1 & 1 \\
1 & 2 & 2  \\
1 & 2 & 2  
\emat
\]
both have limits with coordinates 
\begin {align*} 
a & = \frac{5-\sqrt{17}}{2} = 0.4384471870 \\ 
b& = \frac{-3+\sqrt{17}}{4} =0.2807764065 \\
c & =  \frac{7 - \sqrt{17}}{8} = 0.3596117968.
\end {align*} 
Moreover,  
\[
\lim_{K\rightarrow \infty} S(A_1) = 
\bmat 
1 & 0 & 0 \\
0 & 1/2 & 1/2 \\
0 & 1/2 & 1/2 
\emat.
\]

The matrix 
\[
A_4 = \bmat
1 & K & K \\
K & 1 & 1 \\
K & 1 & 1
\emat
\]
is an $MBN$ matrix with $k=1$, $\ell = 2$, $M=N=1$, and $B = K$.
We have $MN/B^2 = 1/K^2 \neq 0$, and 
\[
A_4 \rightarrow \cdots \rightarrow S(A_4) = 
\bmat 
a & b & b \\
b & c & c \\
b & c & c
\emat
\]
with 
\begin {align*}
a & =  \frac{-K^2 -2+ K \sqrt{K^2 + 8}}{2(K^2-1)} \\
b& =       \frac{3K^2 - K \sqrt{K^2 + 8}}{4(K^2-1)}      \\
c & =      \frac{K^2 - 4 + K \sqrt{K^2 + 8}}{8(K^2-1)}. 
\end {align*}

For example, with $K=2$, we have 
\[
a  = -1 + \frac{2  \sqrt{3}}{3}, \qquad 
b = 1 -  \frac{\sqrt{3} }{3}, \qquad 
c =   \frac{ \sqrt{3} }{6}
\]

Moreover,
\[
\lim_{K\rightarrow \infty} S(A_4) = 
\bmat 
0 & 1/2 & 1/2 \\
1/2 & 1/4 & 1/4 \\
1/2 & 1/4 & 1/4 
\emat.
\]


\section{The matrix $A_5$}
The construction of the Sinkhorn limit of the $3 \times 3$ matrix
\[
A_5 = \bmat 
K & 1 & 1 \\
1 & K & 1 \\
1 & 1 & 1
\emat 
\]
requires only high school algebra.  
There exists a unique positive diagonal matrix $X = \diag(x,y,z)$  
such that $XA_5X$ is doubly stochastic.  We have
\[
S(A_5) = XA_5X = \bmat 
Kx^2 & xy & xz \\
xy & Ky^2 & yz \\
xz & yz & z^2 
\emat
\]
and so 
\begin {align*}
Kx^2 + xy +  xz  & = 1 \\
xy + Ky^2 + yz  & = 1 \\
xz + yz + z^2   & = 1 
\end {align*}
We have 
\[
z = \frac{1 - Kx^2 - xy }{x} = \frac{1 - xy  - Ky^2 }{y}.  
\]
Rearranging, we obtain 
\beq                         \label{Sinkhorn:zC} 
(y-x)((K-1)xy+1) = 0.
\eeq
Note that $0 < xy < 1$. 
If $K > 1$, then  $(K-1)xy+1 > 1$.
 If $0 < K < 1$, then 
\[
0 < (1-K)xy < 1-K < 1
\]
and  $(K-1)xy+1 > 0$.
Therefore, $x=y$, and so 
\beq                         \label{Sinkhorn:zA} 
(K+1)x^2 + xz  = 1 
\eeq
\beq                      \label{Sinkhorn:zB} 
2xz + z^2  = 1.   
\eeq
We obtain   
\[
2\left( 1 - (K+1)x^2 \right)  +  \left(\frac{  1 - (K+1)x^2}{x} \right)^2 = 1.   
\]
Equivalently, 
\[
(K^2-1)x^4 - (2K+1)x^2 + 1 = 0
\]
and so 
\[
x^2 = \frac{2K+1 \pm \sqrt{4K+5}}{2(K^2 -1)}.
\]
Eliminating $xz$ from~\eqref{Sinkhorn:zA} and~\eqref{Sinkhorn:zB} gives 
\[
z^2 = 2(K+1)x^2 - 1 = \frac{K+2 \pm \sqrt{4K+5}}{K-1}.
\]
The inequalities $Kx^2 < 1$ and $z^2 < 1$ imply 
\[
x^2 = \frac{2K+1 - \sqrt{4K+5}}{2(K^2 -1)} 
\]
and 
\[
z^2  = \frac{K+2 - \sqrt{4K+5}}{K-1}.
\]
Thus,
\[
S(A_5) = \bmat
a & b & c \\
b & a & c \\
c & c & d
\emat
\]
where 
\begin {align*}
a & = Kx^2 = \frac{K(2K+1 - \sqrt{4K+5})}{2(K^2 -1)}  \\
b & = x^2  = \frac{2K+1 - \sqrt{4K+5}}{2(K^2 -1)}  \\
c & = xz = \sqrt{ \left(  \frac{2K+1 - \sqrt{4K+5}}{2(K^2 -1)}  \right) \left(  \frac{K+2 - \sqrt{4K+5}}{K-1} \right)   }  \\
d & =   z^2  = \frac{K+2 - \sqrt{4K+5}}{K-1}.
\end {align*}

For example, with $K=2$, we obtain 
\begin {align*}
a & =   \frac{ 5 - \sqrt{13}}{3}  = 0.464816242\\
b & =   \frac{ 5 - \sqrt{13}}{6} =  0.2324081208 \\
c & = \sqrt{bd} = .3027756379 \\
d & =  z^2  = 4- \sqrt{13} = 0.394448725.
\end {align*}

We have the asymptotic limit
\[
\lim_{K\rightarrow \infty} S(A_5) = \bmat 1 & 0 & 0 \\ 0 & 1 & 0 \\ 0 & 0 & 1 \emat.
\]

\section{The matrix $A_6$}

The construction of the Sinkhorn limit of the $3 \times 3$ matrix
\beq          \label{Sinkhorn:KK}
A_6 = \bmat
K & K & 1\\
K & 1 & 1 \\
1 & 1 & 1 
\emat
\eeq
also requires only high school algebra.  
There exists a unique positive diagonal matrix 
$X  = \diag(x,y,z)$ such that 
\[
S(A_6)  = X A_6 X  =  \bmat
Kx^2 & Kx y & x z  \\
Kx y & y^2 & y z  \\
x z & y z & z^2
\emat 
\]
is a doubly stochastic matrix, and so 
\begin{align}
Kx^2 + Kx y + x z  & =1     \label{Sinkhorn:A6-eqn1}              \\
Kx y + y^2 + y z  & = 1  \label{Sinkhorn:A6-eqn2}             \\
x z + y z + z^2 & = 1.  \label{Sinkhorn:A6-eqn3}            
\end{align}
From~\eqref{Sinkhorn:A6-eqn1}, we obtain 
\beq                \label{Sinkhorn:A6-eqn4}      
z = \frac{1}{x} -Kx - Ky.
\eeq
Inserting~\eqref{Sinkhorn:A6-eqn4}  into~\eqref{Sinkhorn:A6-eqn2} gives 
\beq               \label{Sinkhorn:A6-eqn5}      
x = \frac{y}{(K-1)y^2  + 1}.
\eeq
Inserting~\eqref{Sinkhorn:A6-eqn5}  into~\eqref{Sinkhorn:A6-eqn4} gives 
\beq               \label{Sinkhorn:A6-eqn6}      
z = \frac{1}{y} - y - \frac{Ky}{(K-1)y^2 + 1} 
 = \frac{- (K-1)y^4 - 2y^2 + 1}{y( (K-1)y^2+1)}.
\eeq
Inserting~\eqref{Sinkhorn:A6-eqn5} and~\eqref{Sinkhorn:A6-eqn6}  
into~\eqref{Sinkhorn:A6-eqn3} and rearranging  gives 
\[
(K-1)^2 y^6 +3(K-1)y^4 + (K-1)y^2 =1.
\]
Equivalently, 
\[
\left( (K-1)y^2 + 1 \right)^3 = K
\]
and so 
\[
y^2 = \frac{K^{1/3} -1}{K-1} = \frac{1}{1+K^{1/3} + K^{2/3}}
\]
and
\[
y = \frac{1}{\sqrt{1+K^{1/3} + K^{2/3}}}.
\]
Inserting this  into~\eqref{Sinkhorn:A6-eqn5} gives 
\[
x = \frac{y}{K^{1/3}} = \frac{1}{ K^{1/3}\sqrt{1+K^{1/3} + K^{2/3}}}.
\]
and then~\eqref{Sinkhorn:A6-eqn4} gives 
\[
z =  \frac{K^{1/3}}{ \sqrt{1+K^{1/3} + K^{2/3}}}.
\]
Thus,
\[
x^2 = \frac{1}{ K^{2/3}(1+K^{1/3} + K^{2/3})} =  \frac{K^{1/3} -1}{ K^{2/3}(K-1)} 
\]
and
\[
z^2 = \frac{K^{2/3}}{ 1+K^{1/3} + K^{2/3}} =  \frac{ K -K^{2/3}}{ K-1}.
\]
This determines the scaling matrix X.  
The Sinkhorn limit is the circulant matrix 
\[
S(A_6) = \bmat a & b & c \\ b & c & a \\ c & a & b \emat
\]
with 
\begin {align*}
a  & = Kx^2 = yz = \frac{K^{2/3} -K^{1/3} }{K-1} \\
b & = z^2 = Kxy  = \frac{K - K^{2/3} }{K-1} \\
c & = xz = y^2  = \frac{K^{1/3} -1}{K-1}.
\end {align*}
The asymptotic limit is
\[
\lim_{K\rightarrow \infty} S(A_6) = \bmat 0 & 1 & 0 \\ 1 & 0 & 0 \\ 0 & 0 & 1 \emat.
\]

Let 
\[
A_6^{(\ell)} = \bmat a_{i,j}^{(\ell)} \emat
\]
be the $\ell$th matrix in the alternate minimization algorithm for the matrix~\eqref{Sinkhorn:KK}.
We have 
\[
\lim_{\ell\rightarrow \infty} \frac{a_{1,1}^{(\ell)} }{a_{1,3}^{(\ell)} } 
\lim_{\ell\rightarrow \infty} (K-1)a_{1,3}^{(\ell)} +1 
= K^{1/3}
\]
and so alternate minimization generates sequences of rational numbers that converges to $K^{1/3}$.  

For example, 
with $K=2$, we obtain
\[
S(A_6) = \bmat 2^{2/3} - 2^{1/3} & 2 - 2^{2/3}  & 2^{1/3} - 1 \\ 2 - 2^{2/3}  & 2^{1/3} - 1 & 2^{2/3} - 2^{1/3} \\ 2^{1/3} - 1 & 2^{2/3} - 2^{1/3} & 2 - 2^{2/3}  \emat.  
\]


\section{The matrix $A_7$}

Consider the symmetric $3\times 3$ matrix 
\[
A_7= \bmat
K & K & 1 \\
K & 1 & 1 \\
1 & 1 & K
\emat.
\]
There exists a unique positive diagonal matrix $X = \diag(x,y,z)$ such that 
\[
S(A_7) = XA_7X = \bmat
 K x^2 &  K xy & xz \\
 Kxy & y^2 & yz \\
xz & yz & Kz^2
\emat 
\]
is doubly stochastic.  Therefore, 
\begin{align}
 K x^2 +  K xy+ xz  & = 1          \label{Sinkhorn:A7-eqn1}       \\
 Kxy + y^2 + yz  & = 1                 \label{Sinkhorn:A7-eqn2}   \\
xz + yz +  K z^2 & = 1           \label{Sinkhorn:A7-eqn3}  
\end{align}
Observe that equations~\eqref{Sinkhorn:A7-eqn1} and~\eqref{Sinkhorn:A6-eqn1} 
are identical, and that equations~\eqref{Sinkhorn:A7-eqn2} and~\eqref{Sinkhorn:A6-eqn2} 
are identical.  Therefore,
\beq                               \label{SinkhornExplicit:x}
x = \frac{y}{(K-1)y^2 + 1}.
\eeq
and 
\beq                               \label{SinkhornExplicit:z}
z = \frac{- (K-1)y^4 - 2y^2 + 1}{y( (K-1)y^2+1)}.
\eeq
Substituting~\eqref{SinkhornExplicit:x}  and~\eqref{SinkhornExplicit:z} 
into the third equation gives a polynomial in one variable:
\[
(K-1)^3y^8+3(K-1)^2y^6-(K-1)(2K-3)y^4-(4K-1)y^2+K = 0.
\]

By Sinkhorn's theorem, this polynomial has at least one positive solution.
If $K>1$, then, by Descartes's rule of signs, this polynomial has exactly two positive solutions.  
If $0 < K < 1$, then this polynomial has two, four, or six positive solutions.

For example, let $K=2$.
Let $X = \diag(x,y,z)$ be the unique positive diagonal matrix such that the matrix 
\[
S(A_7) = XA_7X = \bmat
2 x^2 & 2 xy& xz \\
2xy & y^2 & yz \\
xz & yz & 2 z^2
\emat 
\]
is doubly stochastic, and 
\begin {align*}
2 x^2 + 2 xy+ xz  & = 1\\
2xy + y^2 + yz  & = 1 \\
xz + yz + 2 z^2 & = 1 
\end {align*}
The number $y$ is a solution of the octic polynomial 
\[
y^8 + 3y^6 - y^4 -7y^2 + 2 = 0.
\]
According to Maple, the unique solution of this polynomial in the interval $(0,1)$ is 
\[
y = 0.533828905923539.
\]
From~equations~\eqref{SinkhornExplicit:x} and~~\eqref{SinkhornExplicit:z}, we obtain 
\[
x = 0.415439687028039
\]
and 
\[
z = 0.508551090023910.
\]
We obtain 
\begin {align*}
a = 2x^2 & =   0.345180267115910   \\
b = 2xy & =    0.443547427206792   \\
c = xz & = 0.211272305677301     \\
d = y^2 & =    0.284973300799523   \\
e = yz  & =  0.271479271993687     \\
f = 2 z^2 & =    0.517248422329014. 
\end {align*}
This agrees with the calculation in Section~\ref{Sinkhorn:section:calculation}.

Let $K= 3$.
Let $X = \diag(x,y,z)$ be the unique positive diagonal matrix such that the matrix 
\[
S(A) = X AX = \bmat
3 x^2 & 3 xy& xz \\
3xy & y^2 & yz \\
xz & yz & 3 z^2
\emat 
\]
is doubly stochastic, and 
\begin {align*}
3 x^2 + 3 xy+ xz  & = 1\\
3xy + y^2 + yz  & = 1 \\
xz + yz + 3z^2 & = 1 
\end {align*}
The number $y$ is a solution of the octic polynomial 
\beq         \label{Sinkhorn:poly-3}       
8y^8+12y^6-6y^4-11y^2+3 = 0.
\eeq
According to Maple, the  solutions of this polynomial in the interval $(0,1)$ are 
\[
0.5083028225 \qand 0.9007108688.
\]
Choosing $y = 0.5083028225$, we obtain from~equations~\eqref{SinkhornExplicit:x} and~\eqref{SinkhornExplicit:z} the numbers 
\[
x = 0.335127736635918
\]
and 
\[
z = 0.453645164346447
\]
and so  
\begin {align*} 
a  = 3x^2 & = 0.336931799588139 \\
b  = 3xy & = 0.511039123248612 \\
 c = xz  & =  0.152029077163254 \\
d = y^2 & = 0.258371759319391\\
e  =  yz & = 0.230589117432000 \\
f  = 3z^2 & =  0.617381805404745
\end {align*}
This agrees with the calculation in Section~\ref{Sinkhorn:section:calculation}.

It is interesting to observe that if we choose the the second root 
of the polynomial~\eqref{SinkhornExplicit:x}, we obtain 
\begin {align*}
x & = .343447174245447 \\
y &:= .900710868780307 \\
z &= -.820818203542269
\end {align*}
and
\begin {align*}
a & =   0.353867884491546 \\
b & =      0.928039808084274 \\
c & =             -0.281907692575816 \\
 d & =                0.811280069138975 \\
e & =            -0.739319877223248 \\
f & =            2.02122756979907
\end {align*}

For matrices of the form $A_7$, we do not explicit formulae for the coordinates 
of the Sinkhorn limit as explict functions of $K$.  Computer calculations suggest 
that the asymptotic limit of $S(A_6)$ as $K\rightarrow \infty$ is 
\[
 \bmat 0 & 1 & 0 \\ 1 & 0 & 0 \\ 0 & 0 & 1 \emat.
\]


\section{Gr\" obner bases and algebraic numbers}
I like solving problems using high school algebra.  
However, it is important to note that the previous calculations are also easily done using Gr\" obner bases. 
 
Here is an example. 
Consider the $A_7$ matrix 
\[
\bmat
K & K & 1 \\
K & 1 & 1 \\
1 & 1 & K
\emat
\]
with $K>0$ and $K \neq 1$.  
There exist unique positive real numbers $x,y,z$ that satisfy the polynomial 
equations 
\begin {align*}
 K  x^2 +  K  xy+ xz  & = 1\\
 K xy + y^2 + yz  & = 1 \\
xz + yz +  K  z^2 & = 1. 
\end {align*}
Equivalently, $(x,y,z)$ is the unique positive vector in $\R^3$ that is in the affine variety 
$V(I)$, where $I$ is the ideal in $\R[x,y,z]$ generated by the polynomials 
\begin {align*}
 K  x^2 +  K  xy+ xz  & - 1\\
 K xy + y^2 + yz  & - 1 \\
xz + yz +  K  z^2 & - 1. 
\end {align*}
Let $K=2$.  
Using the Groebner package in Maple with the lexicographical order $(x,y,z)$, we obtain 
the Gr\" obner basis 
\begin {align*}
f_1(z) & = 4-28 z^2+62 z^4-57 z^6+18 z^8 \\
f_2(y,z) & =  -17 z^3+39 z^5-18 z^7+2 y \\
f_3(x,z) & =   -20 z+96 z^3-135 z^5+54 z^7+4 x 
\end {align*}  
Applying Maple with the lexicographical order $(y,z,x)$, we obtain 
the Gr\" obner basis 
\begin {align*}
g_1(x) & = 2-17 x^2+22 x^4+48 x^6+36 x^8 \\
g_2(x,z) & = -103 x+378 x^3+624 x^5+396 x^7+14 z \\
g_3(x,y) & =  3 x-56 x^3-72 x^5-36 x^7+7 y
\end {align*}  
Applying Maple with the lexicographical order $(z,x,y)$, we obtain 
the Gr\" obner basis 
\begin {align*}
h_1(y) & = 2-7 y^2-y^4+3 y^6+y^8 \\
h_2(x,y) & =  -4 y+2 y^5-3 y^3+y^7+6 x \\
h_3(y,z) & =   -7 y+5 y^5+3 y^3+y^7+6 z
\end {align*}  
Thus, $x^2$, $y^2$, and $z^2$ are algebraic numbers of degree at most 4, and we have explicit polynomial representations of each variable $x$, $y$, $z$ in terms of the others.

 For arbitrary $K$, applying Maple with the lexicographical order $(y,z,x)$, we obtain 
the Gr\" obner basis 
\begin {align*}
h_1(y) 
& = K- (4 K - 1) y^2 -(K-1)(2K-3) y^4+3 (K-1)^2  y^6+(K-1)^3 y^8 \\
h_2(x,y) 
& =   K (K+1) x  -2 K y - (K-1)(2K-1) y^3+ 2 (K-1)^2  y^5+(K-1)^3 y^7 \\
h_3(y,z) 
& = K (K+1) z - (K-1)^2  y -3(K-1) y^3+ (K-1)^2(K-3) y^5+( K-1)^3 y^7.
\end{align*}   
For each of the 8 roots of $h_1(y)$,the polynomials $g_2(z,y)$ and $g_3(x,y)$ determine unique numbers $x$ and $z$.  
Exactly one of the triples $(x,y,z)$ will be positive.  

For every positive symmetric $n \times n$ matrix $A = (a_{i,j})$, 
the Sinkhorn limit  $S(A) = XAX$ with scaling matrix $X = \diag(x_1,\ldots, x_n)$ 
is the unique positive solution of a set 
$Q = \{q_i: i=1,\ldots, n\}$  of $n$ quadratic equations of the form 
\[
q_i = q_i(x_1,\ldots, x_n) = \sum_{j=1}^n a_{i,j} x_ix_j - 1 = 0.
\]
 Equivalently, $(x_1,\ldots, x_n)$ is the unique positive vector in the affine variety of the ideal generated by $Q$.
A Gr\" obner basis for this ideal shows that  if the coordinates of the matrix $A = (a_{i,j})$ 
are rational numbers, then $x_1,\ldots, x_n$ are 
algebraic numbers of degrees bounded in terms of $n$.


\section{Diophantine approximation}

Let $A$ be a an $n \times n$  matrix with positive rational coordinates, and let $d$ be the least 
common multiple of the denominators of the coordinates of $ A $.  
The matrix $dA $ has positive integral coordinates, and  the matrix obtained 
by row scaling (or column scaling) $ A $ is equal to the matrix obtained 
by row scaling (or column scaling) $dA $.  Thus, the Sinkhorn limit obtained from the rational matrix $ A $ 
equals  the Sinkhorn limit obtained from the integral matrix $d A $.   
The sequence of matrices generated by alternate row and column scalings are rational matrices.  
If $A^{(\ell)} = \left( a_{i,j}^{(\ell)} \right)$ is the $\ell$th matrix obtained in the alternate minimization 
algorithm, and if the Sinkhorn limit is $S(A) = \left( s_{i,j} \right) $, then 
\[
\lim_{\ell\rightarrow \infty} a_{i,j}^{(\ell)} = s_{i,j}
\]
for all $i,j = 1,\ldots n$.  If the coordinate $s_{i,j}$ is irrational for some pair $(i,j)$, 
then the alternate minimization cannot terminate in a finite number of steps.  
It is an open problem to the matrices $A$ for which the alternate minimization does 
terminate in a finite number of steps.  

The Sinkhorn limit coordinates $s_{i,j}$ are algebraic numbers for all rational matrices A.  
If the coordinate $s_{i,j}$ is irrational for some $i$ and $j$, then the alternate minimization 
algorithm constructs a sequence of rational approximations to $s_{i,j}$.  
For example, alternate minimization provides a sequence (in fact, several sequences) 
of rational numbers that converge to $K^{1/3}$  for every positive integer $K$.  
The  matrix 
\[
A_6 = \bmat
K & K & 1\\
K & 1 & 1 \\
1 & 1 & 1 
\emat
\]
has Sinkhorn limit  
\[
S(A_6) = \bmat a & b & c \\ b & c & a \\ c & a & b \emat
\]
with 
\begin {align*}
a  & =  \frac{K^{2/3} -K^{1/3} }{K-1} \\
b & =  \frac{K - K^{2/3} }{K-1} \\
c & = \frac{K^{1/3} -1}{K-1}.
\end {align*}
If $A_6^{\ell)} = \left( a_{i,j}^{\ell)} \right)$, then 
\begin{align*}
K^{1/3} -1
& = \lim_{\ell \rightarrow \infty} (K-1)  a_{1,3}^{(\ell)}  \\
& = \lim_{\ell \rightarrow \infty} (K-1)  a_{2,2}^{(\ell)}  \\
& = \lim_{\ell \rightarrow \infty} (K-1)  a_{3,1}^{(\ell)}.  
\end{align*} 
For example, for $K=2$, we have 
\[
2^{1/3} -1
= \lim_{\ell \rightarrow \infty}  a_{1,3}^{(\ell)}  \\
= \lim_{\ell \rightarrow \infty}  a_{2,2}^{(\ell)}  \\
= \lim_{\ell \rightarrow \infty} a_{3,1}^{(\ell)}.  
\]
Here are the rational numbers in the first six iterations of the Sinkhorn algorithm, 
and their decimal representations:

\begin{center}
\begin{tabular}{c | c | c | c}  
$ \ell$ & $a^{( \ell)}_{1,3}$ & $a^{( \ell)}_{2,2}$  & $a^{( \ell)}_{3,1}$  \\ \hline 
&&&\\
$1$ & $ \frac{ 1 }{ 5 }  $ & $   \frac{  1}{ 4 }    $  & $    \frac{ 1 }{  3}     $  \\
&&&\\
$2$ & $ \frac{ 12 }{ 47 }  $ & $   \frac{ 15 }{  59}    $  & $    \frac{10  }{ 37 }     $  \\  
&&&\\
$3$ & $ \frac{ 2183 }{ 8434 }  $ & $   \frac{ 1739 }{ 6695}    $  & $    \frac{2773}{ 10617}     $  \\   
&&&\\
$4$ & $ \frac{71080815}{273555853} $ & $   \frac{ 44771889}{172318334 }    $  & $    \frac{ 56465630}{217090223}     $  \\
&&&\\
$5$ & $   \frac{ 37408625555048482}{143933615530682603 }    $  & $    \frac{ 59386301130725219}{228480929930639987}     $ 
& $ \frac{ 47138688844908902} {181342241085731085 }  $  \\
&&&\\
$6$ & $ \frac{ 41433243878974147831553607829895895}{159407905344245227309688080035616727 }  $ & $   \frac{ 26101244407905972515151593345814255}{100420574611609687570620843932756311  }    $  & 
$  a^{(6)}_{3,1}  $  \\&&&\\
\end{tabular}
\end{center}
where
\begin{align*}
a^{(6)}_{3,1} 
& =   \frac{ 32886086324729567223915642757046161}{126521819019515660085772437278570566}   \\ 
& \\
& = 0.2599242295\ldots
\end{align*}
Note that 
\[
2^{1/3}-1 = 0.2599210499\ldots.
\]
The continued fraction for $2^{1/3}-1$ is
$[0, 3, 1, 5, 1, 1, 4, 1, 1, 8, 1, 14, 1, 10,\ldots] .$
For comparison, here are the first ten convergents of the continued fraction for $2^{1/3}-1$:
\[
\begin{matrix}
 \frac{1}{3}  & =              0.3333333333  & \hspace{0.8cm} &   
 \frac{1}{4}  & =              0.2500000000  \\
 &&&\\
 \frac{6}{23}  & =            0.2608695652  & \hspace{0.8cm} &   
 \frac{7}{27}    & =          0.2592592593 \\
 &&&\\
\frac{13}{50} & =          0.2600000000  & \hspace{0.8cm} &  
\frac{59}{227} & =        0.2599118943 \\
&&&\\
\frac{72}{277} & =       0.2599277978  & \hspace{0.8cm} &   
\frac{131}{504} & =     0.2599206349 \\
&&&\\
\frac{1120}{4309} & = 0.2599210954  & \hspace{0.8cm} &    
\frac{1251}{4813} & = 0.2599210472
\end{matrix}
\]


\section{Rationality and finite length}
For what positive $n\times n$ matrices does the alternate minimization algorithm 
converge in finitely many steps?   This problem has been solved for $2 \times 2$ matrices 
(Nathanson~\cite{nath18aaa}), but it is open for all dimensions $n \geq 3$.  
In dimension 3, matrices equivalent to $A_1$ become doubly stochastic in one step, 
that is, after one row or one column scaling.  
It is not know if there exists a positive  $3\times 3$ matrix that becomes doubly stochastic in exactly two steps.
More generally, it is not know if there exists  a positive  $3\times 3$ matrix 
that becomes doubly stochastic in exactly $s$ steps for some $s \geq 2$.

Consider the matrix $A_2 = \bmat K & 1 & 1 \\ 1 & 1 & 1 \\ 1 & 1 & 1 \emat$ with parameter $K$.   
If $K$ is a rational number, then every matrix generated by iterated row 
and column scalings has rational coordinates.  
If the Sinkhorn limit contains an irrational coordinate, 
then the alternate minimization algorithm cannot terminate in finitely many steps.  

If $K$ is an integer and $K \geq 2$, 
then the Sinkhorn limit $S(A_2)$ has coordinates in the quadratic  field 
 $\Q(\sqrt{8K+1})$.  For example, from~\eqref{Sinkhorn:A2-a},  
 the $(1,1)$ coordinate of $S(A_2)$ is 
 \[
 \frac{2K+1-\sqrt{8K+1}}{2(K-1)}.
 \]
This number is rational if and only if the odd integer $8K+1$ is the  square of an odd integer, that is,
if and only if $8K+1 = (2r+1)^2$ for some positive integer $r$ 
and  so $K = r(r+1)/2$ is a triangular number.  
From~\eqref{Sinkhorn:A2-a},~\eqref{Sinkhorn:A2-b}, and~\eqref{Sinkhorn:A2-c}, 
we obtain 
\begin {align*}
a & = \frac{r^2-r}{r^2+r-2} = \frac{r}{r+2}  \\
b & = \frac{r -1}{r^2+r-2}  = \frac{1}{r+2} \\
c & =   \frac{r^2-1}{2(r^2+r-2)} = \frac{r+1}{2(r+2) }.
\end {align*}
Moreover, $S(A_2) = X A_2 X$, where $X = \diag (x,y,y)$ 
with $Kx^2 = a$ and $y^2 = c$.  Thus, 
\[
x = \sqrt{\frac{a}{K} } =  \sqrt{ \frac{2}{ (r+1)(r+2)}  }
\qqand 
y = \sqrt{c} =  \sqrt{   \frac{r+1}{2(r+2)}    }.  
\]

For example, if $K = 3$, then $r=2$ and 
\[
A_2 = \bmat 3 & 1 & 1 \\ 1 & 1 & 1  \\ 1 & 1 & 1 \emat
\rightarrow \cdots \rightarrow 
X A_2 X = S(A_2 ) = \bmat 1/2 & 1/4 & 1/4 \\ 1/4 & 3/8 & 3/8  \\ 1/4 & 3/8 & 3/8 \emat 
\]
where 
\[
 X = \diag( \sqrt{6}/6, \sqrt{6}/4,  \sqrt{6}/4 ).
\]
Note that  $A_2$ also has a scaling by rational matrices  
\[
S(A_2) = X'A_2 Y'
\]
where 
\[
X' = \diag( 1/6, 1/4,1/4) \qqand Y' = \diag(1, 3/2, 3/2).
\]
It is not known if there exists a triangular number $K$ 
for which the alternate minimization  algorithm terminates in a finite number of steps.

\section{Open problems}

\benum

\item
Compute explicit formulas for the Sinkhorn limits of all positive symmetric $3\times 3$ matrices.
This is a central problem.

\item
Here is a special case.  
Let $K, L, M$ and 1 be pairwise distinct positive numbers.  Compute the Sinkhorn limits of the matrices
\[
\bmat 
K & 1 & 1 \\
1 & L & 1 \\
1 & 1 & 1
\emat
\qqand
\bmat 
K & 1 & 1 \\
1 & L & 1 \\
1 & 1 & M
\emat.
\]

\item
For what positive $n\times n$ matrices does the alternate minimization algorithm 
converge in finitely many steps?   This is the problem discussed in the previous section.

\item
It is not known what algebraic numbers appear as coordinates of the Sinkhorn limit 
of a positive integral matrix.  
It would be interesting to have an example of an algebraic number in the unit interval 
that is not a coordinate of the Sinkhorn limit of a rational matrix.

\item
Does there exist a $3\times 3$ matrix $A$ such that $A$ is row stochastic but not column stochastic, 
and $AY(A)$ is doubly stochastic?

\item
Does every possible shape  of a doubly stochastic $3\times 3$ matrix $A$ appear as the 
nontrivial limit of some $3\times  3$ matrix?

\item
Why does the shape of the Sinkhorn limit $S(A)$ seem to depend only on the shape of the matrix $A$ 
and not on the numerical values of the coordinates of $A$?
\item
What does the Sinkhorn limit  $S(A)$ tell us about the matrix $A$?  
What information does it convey? 

\item
The matrix $ A $ is \emph{positive} if $a_{i,j} > 0$ for all  $i$ and  $j$.
The matrix $ A $ is \emph{nonnegative} if $a_{i,j} \geq 0$ for all  $i$ and  $j$.

Let A\ be a nonnegative $m\times n$ matrix.
Let $ \mbr = (r_1, r_2,\ldots, r_m) \in \R^m$ and let $ \mbc = (c_1, c_2, \ldots, c_n) \in \R^n$.
The matrix A\ is \emph{$ \mbr$-row stochastic} 
if $\row_i( A ) = r_i$ for all $i \in \{1,2,\ldots, m\}$.
The matrix A\ is \emph{$ \mbc$-column stochastic} 
if $\col_j( A ) =c_j$ for all $j \in \{1,2,\ldots, n\}$.
The matrix $ A $ is \emph{$(\mbr,\mbc)$-stochastic}  
if it is both $\mbr$-row stochastic and $\mbc$-column stochastic.
Note that if A\ is $(\mbr,\mbc)$-stochastic, then 
\beq           \label{Sinkhorn:rc-condition}
\sum_{i=1}^m r_i = \sum_{i=1}^m \sum_{j=1}^n a_{i,j} 
=  \sum_{j=1}^n  \sum_{i=1}^ma_{i,j} = \sum_{j=1}^n c_j.
\eeq
Let A\ be a positive matrix.  Let $X$ be the $m\times m$ diagonal matrix whose 
$i$th coordinate is $r_i/\row_i(A)$, and let $Y$ be the $n \times  $ diagonal matrix whose 
$j$th coordinate is $c_j/\col_j(A)$.
The matrix $XA$ is $\mbr$-row stochastic and the matrix 
$AY$ is $\mbc$-column stochastic.  

A simple modification of the alternate minimization algorithm applied to a positive matrix 
satisfying~\eqref{Sinkhorn:rc-condition} produces 
an $(\mbr,\mbc)$-stochastic Sinkhorn limit.  
It is an open problem to compute explicit Sinkhorn limits in 
the $(\mbr,\mbc)$-stochastic setting.

\eenum

\section{Notes}

In his 1964 paper, Richard Sinkhorn~\cite[p.877]{sink64} wrote:
\begin{quotation}
The iterative process of alternately normalizing the rows and columns of a
strictly positive $N \times N$ matrix is convergent to a strictly positive 
doubly stochastic matrix.
\end{quotation}
Sinkhorn did not prove this result.  
The proof of convergence of the alternate minimization algorithm appears 
in  Knopp and Sinkhorn~\cite{sink-knop67}, 
and in Letac~\cite{leta74}.
Geometric existence proofs of exact scaling appear in Menon~\cite{meno67}, 
and  in Tverberg~\cite{tver76}.  

The computational complexity of Sinkhorn's alternate scaling algorithm 
is investigated in 
Kalantari and Khachiyan~\cite{kala-khac93,kala-khac96},
Kalantari,  Lari, Ricca, and Simeone~\cite{kala08}, 
Linial, Samorodnitsky and Wigderson~\cite{lini-samo-wigd98b}
and Allen-Zhu, Li, Oliveira, and Wigderson~\cite{alle-li-oliv-wigd17}.
An extension of matrix scaling to operator scaling began 
with Gurvits~\cite{gurv04}, and is developed in 
Garg,  Gurvits, Oliveira, and Wigderson~\cite{garg-gurv-oliv-wigd16,garg-gurv-oliv-wigd17},
Gurvits~\cite{gurv15}, and Gurvits and Samorodnitsky~\cite{gurv-samo14}.  
Motivating some of this recent work are the classical papers of Edmonds~\cite{edmo67} 
and Valient~\cite{vali79a,vali79b}.

The literature on matrix scaling is vast.  See the recent survey paper of 
Idel~\cite{idel16}.  For the early history of matrix scaling, see  
Allen-Zhu, Li, Oliveira, and Wigderson~\cite[Section 1.1]{alle-li-oliv-wigd17}.

\emph{Acknowledgements.}
The alternate minimization algorithm was discussed in several lectures 
in the New York Number Theory Seminar, 
and I thank the participants for their useful remarks.  
In particular, I thank David Newman for making the initial computations 
that suggested some of the problems considered in this paper.

\def\cprime{$'$} \def\cprime{$'$} \def\cprime{$'$}
\providecommand{\bysame}{\leavevmode\hbox to3em{\hrulefill}\thinspace}
\providecommand{\MR}{\relax\ifhmode\unskip\space\fi MR }
\providecommand{\MRhref}[2]{%
  \href{http://www.ams.org/mathscinet-getitem?mr=#1}{#2}
}
\providecommand{\href}[2]{#2}

\end{document}